\documentstyle{amsppt}
\raggedbottom
\NoBlackBoxes
\nologo
\tolerance=1000
\define\Polacik{Pol\'a\v{c}ik\ }
\define\Terescak{Tere\v{s}\v{c}\'ak\ }

\define\JDDE{J. Dynamics Differential Equations}
\define\JDE{J. Differential Equations}
\define\eg{e\.~g\.\ }
\define\resp{resp\.\ }
\define\complex{\Bbb C}
\define\reals{\Bbb R}
\define\integers{\Bbb Z}
\define\naturals{\Bbb N}

\define\wcl{\operatorname{w^{*}-cl}}
\define\codim{\operatorname{codim}}

\define\Cone{C^{1}(\overline\Omega)}

\define\bs{b\cdot s}
\define\bt{b\cdot t}
\define\pshat{\hat\psi}
\define\Pshat{\widehat\Psi}
\define\Bb{\Bbb B}
\define\CalB{\Cal B}
\define\CalBhat{\widehat{\CalB}}
\define\CalBi{\CalB_{i}}
\define\CalBone{\CalB_{1}}
\define\CalBtwo{\CalB_{2}}
\define\CalS{\Cal S}
\define\CalT{\Cal T}
\define\CalZ{\Cal Z}
\define\CalZhat{\widehat{\Cal Z}}
\define\Id{\roman{Id}}

\define\Lone{L^{1}(\Omega)}
\define\Lp{L^{p}(\Omega)}
\define\Lq{L^{q}(\Omega)}
\define\Lr{L^{r}(\Omega)}
\define\Mpbs{{\bold M}_{p}(\bs)}
\define\Mpbt{{\bold M}_{p}(\bt)}
\define\Rn{{\reals}^{n}}
\define\Xst{X^{*}}
\define\0{\{0\}}
\redefine\Re{\operatorname{Re}}
\redefine\ge{\geqslant}
\redefine\le{\leqslant}

\topmatter
\title \nofrills Globally positive solutions of linear
parabolic PDEs of second order with Robin boundary
conditions
\endtitle
\author Janusz Mierczy\'nski
\footnotemark"*"
\endauthor
\footnotetext""{
\copyright\ 1997. This manuscript version is made available under the CC-BY-NC-ND 4.0 license http://creativecommons.org/licenses/by-nc-nd/4.0/}
\footnotetext""{
Published in {\it Journal of Mathematical Analysis and Applications} {\bf 209}(1) (1997), pp. 47--59.}
\footnotetext""{
https://doi.org/10.1006/jmaa.1997.5323}
\footnotetext"*"{
Research supported by KBN grant 2 P03A 076
08.
}
\leftheadtext\nofrills{\eightsmc Janusz Mierczy\'nski}
\rightheadtext\nofrills{\eightsmc Globally positive
solutions of linear parabolic PDEs}
\address
Institute of Mathematics
\newline\indent Technical University of Wroc{\l}aw
\newline\indent Wybrze\D{z}e Wyspia\'nskiego 27
\newline\indent PL-50-370 Wroc{\l}aw
\newline\indent Poland
\endaddress
\email mierczyn\@banach.im.pwr.wroc.pl \endemail
\endtopmatter
\document

\bigskip
Let $\Omega\subset\Rn$ be a bounded domain with boundary
$\partial\Omega$ of class $C^{2}$.

We consider a linear parabolic differential equation of
second order
$$
u_{t}=\sum\limits_{i,j=1}^{n}
\frac{\partial}{\partial x_{i}}
\left(a_{ij}(x)\frac{\partial u}{\partial x_{j}}\right)
+a_{0}(t,x)u,
\qquad t\in\reals, x\in\Omega, \tag{1a}
$$
where $a_{ij}=a_{ji}\in C^{1}(\overline\Omega)$,
$\sum\limits_{i,j=1}^{n}a_{ij}(x)\xi_{i}\xi_{j}>0$ for each
$x\in\overline\Omega$ and each nonzero
$(\xi_{1},\dots,\xi_{n})\in\Rn$, and
$a_{0}\in L^{\infty}(\reals\times\Omega)$.

Equation (1a) is complemented with regular oblique (Robin)
homogeneous boundary conditions
$$
\frac{\partial u}{\partial\nu}(t,x)+c(x)u(t,x)=0, \qquad
t\in\reals, x\in\partial\Omega, \tag{1b}
$$
where $\nu\in C^{1}(\partial\Omega,\Rn)$ is a vector field
pointing out of $\Omega$ and $c\in C^{1}(\partial\Omega)$ is
a nonnegative function.

Motivated by deep results on one-dimensional parabolic
PDE's contained in a series of papers [4], [5] by S.-N.
Chow, K. Lu and J. Mallet-Paret, I investigate the set of
solutions to (1a)+(1b) that are globally positive.

Our main result is the following (see Corollary 2.4).
\proclaim{Theorem}
Let $u_{1}$, $u_{2}$ be solutions of (1a)+(1b) defined and
positive for all $t\in\reals$ and all
$x\in\overline\Omega$.  Then there is a positive constant
$\kappa$ such that $u_{1}(t,x)={\kappa}u_{2}(t,x)$ for all
$t\in\reals$ and all $x\in\overline\Omega$.
\endproclaim

A key idea in the proofs is to consider a family of linear
parabolic PDEs obtained by letting the coefficient
$a_{0}(t,x)$ vary in some metrizable compact subset of the
vector space $L^{\infty}(\reals\times\Omega)$ endowed with
the weak$^{*}$ topology.  That family generates in a
natural way a compact linear skew-product semiflow on a
(product) vector bundle with generic fiber
$L^{1}(\Omega)$.

For $n$ arbitrary, the results obtained seem to be new even
in the case of autonomous equations (compare the remarks on
p\.~288 in [5]).

I would like to mention that it was John Mallet-Paret who
asked me about an ``easy'' way to establish uniqueness (up
to multiplication) of globally positive solutions.  I am
grateful for his hospitality during my short stay at the
Lefschetz Center for Dynamical Systems, Brown University.

\head 0. Preliminaries \endhead

Throughout the paper, for a Banach space $X$ the symbol
${\Cal L}(X)$ stands for the Banach space of bounded linear
maps from $X$ into itself, endowed with the uniform
operator topology.

Let $\Cal A$ be the differential operator
$$
{\Cal A}u(x):=-\sum\limits_{i,j=1}^{n}
\frac{\partial}{\partial x_{i}}
\left(a_{ij}(x)\frac{\partial u}{\partial x_{j}}\right)
\tag{0.1}
$$
Denote by $\tilde{A}_{0}$ the linear operator
$\tilde{A}_{0}u:={\Cal A}u$ acting on the Banach space of
functions from $C^{2}(\overline\Omega)$ satisfying the
boundary condition (1b).  For each $1\le p\le\infty$ the
operator $\tilde{A}_{0}$ is closable and densely defined in
$\Lp$.  Let $\tilde{A}_{p}$ denote the closure of
$\tilde{A}_{0}$ in $\Lp$.  We write $\|\cdot\|_{p}$ for the
norm on the Banach space $\Lp$, $1\le p\le\infty$.

It is well known (see \eg Pazy [10]) that for each $1\le
p<\infty$ the (unbounded) linear operator $\tilde{A}_{p}$
is the infinitesimal generator of an analytic semigroup
$\{e^{-\tilde{A}_{p}t}\}_{t\ge 0}$ of compact linear
operators from ${\Cal L}(\Lp)$.  Put
$A_{p}:=\tilde{A}_{p}+\Id$.

\proclaim{Lemma 0.1}
For each $1\le p<\infty$ the spectrum $\sigma(A_{p})$ of
$A_{p}$ lies in the halfspace
$\{\zeta\in\complex:\Re{\zeta}\ge1\}$.
\endproclaim
\demo{Proof}
As $-A_{p}$ is the infinitesimal generator of an analytic
semigroup $e^{-A_{p}t}$ of compact linear operators, for
each $t>0$ the spectrum $\sigma(e^{-A_{p}t})$ consists of
eigenvalues and $\0$.  By the standard regularity theory
the solutions of the abstract parabolic equation
$u_{t}+A_{p}u=0$ are classical ones, hence we can consider
$e^{-A_{p}\cdot1}\in{\Cal L}(C(\overline\Omega))$.  The
Kre\u{\i}n--Rutman theorem (see \eg Hess [8]) implies that
the spectral radius of $e^{-A_{p}\cdot1}$ is a positive
eigenvalue corresponding to a positive eigenfunction
$v_{0}$.  Multiplying $v_{0}$ by a positive constant, if
necessary, we can assume that $v_{0}(x)\le1$ for each
$x\in\overline\Omega$.

For $x\in\overline\Omega$, $t\ge0$, put
$\tilde{v}(t,x):=e^{-t}$.  The function $\tilde{v}$ is a
supersolution (for the definition see \eg Hess [8] or Smith
[14]) for the initial boundary value problem
$$
\align
&u_{t}+{\Cal A}u=-u,\\
&\frac{\partial u}{\partial\nu}+c(x)u=0,\\
&u(0,x)=1\text{ for }x\in\overline\Omega.\\
\endalign
$$
An application of the parabolic maximum principle to
$\tilde{v}$ and to the solution $v(t,x)$ to the initial
boundary value problem
$$
\align
&v_{t}+{\Cal A}v=-v,\\
&\frac{\partial v}{\partial\nu}+c(x)v=0,\\
&v(0,x)=v_{0}(x)\text{ for }x\in\overline\Omega.\\
\endalign
$$
yields that $v(t,x)$ converges to zero as $t\to\infty$,
uniformly in $x\in\overline\Omega$, with the exponential
rate of decay not larger than $-1$.  So the spectral radius
of $e^{-A_{p}\cdot1}\in{\Cal L}(C(\overline\Omega))$ (hence
of $e^{-A_{p}\cdot1}\in{\Cal L}(\Lp)$) is not larger than
$e^{-1}$. As the spectral mapping theorem (see \eg [10])
asserts that
$\sigma(e^{-A_{p}\cdot1})=\exp(\sigma(-A_{p}))\cup\0$, the
desired result follows.
\qed\enddemo

For $\alpha\ge 0$ denote by $A_{p}^{\alpha}$ the fractional
power of $A_{p}$, and by $\Lp^{\alpha}$ its domain
endowed with the graph norm
$\|\cdot\|_{p,\alpha}:=\|A_{p}^{\alpha}\cdot\|_{p}$.
Obviously, $A_{p}^{1}=A_{p}$ and $\Lp^{1}$ is the domain
of $A_{p}$.

\proclaim{Theorem 0.2}
For each $1\le p<\infty$ the following holds:

{\rm (a)} there is a positive constant $C=C(p)$ such that
$$
\|e^{-A_{p}t}u\|_{p}\le C\|u\|_{p}
\qquad\text{for all }u\in\Lp\text{ and }t\ge0,
$$

{\rm (b)} for each $u\in\Lp$, $t>0$ and $\alpha\ge 0$,
$e^{-A_{p}t}u\in\Lp^{\alpha}$,

{\rm (c)} for each $0\le\alpha\le 1$ there is a positive
constant $C_{1}=C_{1}(p,\alpha)$ such that
$$
\|e^{-A_{p}t}u\|_{p,\alpha}\le
\frac{C_{1}}{t^{\alpha}}\|u\|_{p}
\qquad\text{for all }u\in\Lp\text{ and }t>0,
$$

{\rm (d)} for each $0<\alpha\le 1$ there is a positive
constant $C_{2}=C_{2}(p,\alpha)$ such that
$$
\|(\Id-e^{-A_{p}t})u\|_{p}\le
C_{2}t^{\alpha}\|u\|_{p,\alpha}
\qquad\text{for all }u\in\Lp^{\alpha}\text{ and }t>0.
$$
\endproclaim
\demo{Proof}
All the facts stated above are standard (see \eg Henry [7]
or Pazy [10])
\qed\enddemo

For $t\in\reals$ and $b\in L^{\infty}(\reals\times\Omega)$,
denote by $\bt$ the translate of $b$,
$(\bt)(s,x):=b(t+s,x)$.  Put
$$
\Bb:=\wcl\{a_{0}\cdot t:t\in\reals\},
$$
where $\wcl$ stands for the closure in the weak$^{*}$
topology of $L^{\infty}(\reals\times\Omega)$.  Since $\Bb$
is bounded, it is a compact
metrizable space in the weak$^{*}$ topology of
$L^{\infty}(\reals\times\Omega)$.

For $b\in\Bb$ denote by $\Mpbt$ the multiplication operator by
$b(t,\cdot)$:
$$
(\Mpbt\phi)(x):=b(t,x)\phi(x), \qquad \phi\in\Lp, x\in\Omega.
$$
The following result is straightforward.
\proclaim{Lemma 0.3}
For each $b\in\Bb$, each $1\le p<\infty$ and each
$\phi\in\Lp$ we have $\|\Mpbt\phi\|_{p}\le R\|\phi\|_{p}$,
where $R$ is the essential supremum of $a_{0}$.
\endproclaim

\head 1. Linear skew-product (semi)dynamical systems
\endhead
We consider the following abstract parabolic equation
$$
u_{t}+A_{p}u=\Mpbt u, \tag{1.1}
$$
together with an initial condition
$$
u(0)=u_{0}. \tag{1.2}
$$
In the present section we investigate, for $1\le p<\infty$,
the dependence of solutions to (1.1) satisfying (1.2) on
the coefficient $b\in\Bb$ and the initial condition
$u_{0}\in\Lp$.

The abstract initial value problem (1.1)+(1.2) can be
written in the form of the following integral equation
$$
u(t)=e^{-A_{p}t}u_{0}+
\int\limits_{0}^{t}e^{-A_{p}(t-s)}\Mpbs u(s)\,ds. \tag{1.3}
$$
We refer to solutions of (1.3), that is, functions $u\in
C([0,\infty),\Lp)$ satisfying (1.3) for each $t\ge0$ as
{\it mild solutions\/} of (1.1)+(1.2).  If a (mild) solution
is defined for all $t\in\reals$, we say it is a {\it global
solution\/}.

\proclaim{Proposition 1.1}
For each $b\in L^{\infty}(\reals\times\Omega)$, each $1\le
p<\infty$ and each $u_{0}\in\Lp$ there exists precisely one
solution $u(\cdot;u_{0},b):[0,\infty)\to\Lp$ to (1.3).
\endproclaim
\demo{Proof}
Applying the contraction mapping principle yields the
desired solution.  For a similar result, see Problem~7.1.5
in Henry's book [7].
\qed\enddemo

The proof of the following proposition rests on the
estimates occurring in the proofs of Lemma~3.2 and
Thm\.~3.3 in [5], and we did not repeat it here.

\proclaim{Proposition 1.2}
{\rm a)} For $T>0$ and $1\le p<\infty$ the mapping
$\Lp\times\Bb\ni(u_{0},b)\mapsto u(\cdot;u_{0},b)\in
C([0,T],\Lp)$ is continuous.

{\rm b)} For $0<t_{1}\le t_{2}$, $1\le p<\infty$ and
$0\le\alpha<1$ the mapping $\Lp\times\Bb\ni(u_{0},b)\mapsto
u(\cdot;u_{0},b)\in C([t_{1},t_{2}],(\Lp)^{\alpha})$ is
continuous.
\endproclaim

\proclaim{Theorem 1.3}
{\rm a)} For $t>0$, $u_{0}\in\Lone$ and $b\in\Bb$ we have
$u(t;u_{0},b)\in\Cone$.

{\rm b)} For $0<t_{1}\le t_{2}$ the mapping
$\Lone\times\Bb\ni(u_{0},b)\mapsto u(\cdot;u_{0},b)\in
C([t_{1},t_{2}],\Cone)$ is continuous.
\endproclaim
\demo{Proof}
By Amann [2] or Pazy [10], the domain of $A_{1}$ embeds
continuously in the Sobolev space $W^{1,1}(\Omega)$.  A
corollary of the Nirenberg--Gagliardo inequality (see \eg
[7]) states that
$$
\Lone^{1/2}\subset\Lq\quad\text{if }-n/q\le1/2-n,
$$
that is, $\Lone^{1/2}\subset\Lq$ for each
$1<q<2n/(2n-1)$.  From Proposition~1.2 we deduce that
$u(t;u_{0},b)\in\Lq$ for each $t>0$, each $u_{0}\in\Lone$
and each $b\in\Bb$.

By [10], for $1<p<\infty$ the domain of $A_{p}$ embeds
continuously in the Sobolev space $W^{2,p}(\Omega)$.  A
corollary of the Nirenberg--Gagliardo inequality states
that
$$
\Lp^{1/2}\subset\Lq\quad\text{if }-n/q\le1-n/p.
$$
Therefore, $\Lp^{1/2}\subset\Lq$ for each $p<q<np/(n-p)$.
Consequently, $\Lp^{1/2}\subset\Lr$, where $r:=n/(n-1)$.
Repeating the reasoning sufficiently many times, we prove
that $u(t;u_{0},b)\in\Lp$ for each $t>0$, each
$u_{0}\in\Lone$ and each $b\in\Bb$ and $1\le p<\infty$.

For $p>3n$ a corollary of the Nirenberg--Gagliardo
inequality states that
$$
\Lp^{3/4}\subset C^{\mu}(\overline\Omega)
\quad\text{if }0\le\mu<3/2-1/3=7/6.
$$
It follows that $u(t;u_{0},b)\in\Cone$
for each $t>0$, each $u_{0}\in\Lone$ and each $b\in\Bb$.

Part b) follows similarly by Proposition~1.2.
\qed\enddemo

We recall now the definition of a skew-product dynamical
system (with discrete time).

Let $B\times X$ be a product Banach bundle, where the base
space $B$ is a compact metrizable space and the fiber $X$
is a Banach space.  For a homeomorphism $\phi:B\to B$ the
iterates ${\phi}^{k}$, $k\in\integers$, form a (discrete
time) dynamical system on $B$.  A {\it compact linear
skew-product dynamical system\/}
$\{{\Psi}^{k}\}_{k=1}^{\infty}$ on $B\times X$, {\it
covering\/} $\phi$, is given by a family of compact linear
operators $\{\psi(b):b\in B\}$ depending continuously on
$b\in B$ in the uniform operator topology, in the following
way
$$
\Psi(b,u)=(\phi(b),\psi(b)u), \quad b\in B, u\in X.
$$
In other words, $\Psi$ is a vector bundle endomorphism. The
iterates ${\Psi}^{k}$, $k\in\naturals$, are given by
$$
{\Psi}^{k}(b,u)=({\phi}^{k}(b),{\psi}^{(k)}(b)u),
$$
where we denote
$$
{\psi}^{(k)}(b):= \psi({\phi}^{k-1}(b))\circ
\psi({\phi}^{k-2}(b))\circ\dots\circ
\psi(\phi(b))\circ\psi(b).
$$
(the {\it cocycle identity}).

Let $\Xst$ stand for the Banach space dual to $X$. For a
compact linear skew-product dynamical system
$\{\Psi^{k}\}_{k=1}^{\infty}$ on $B\times X$ we define its
{\it dual system\/} $\{(\Psi^{*})^{k}\}$ on $B\times\Xst$
by
$$
\Psi^{*}(b,u^{*}):=(\phi^{-1}(b),\psi^{*}(b)u^{*}),
\quad b\in B, u^{*}\in\Xst,
$$
where $\psi^{*}(b)\in{\Cal L}(\Xst)$ is the dual operator
to $\psi(b)$.  It is straightforward that
$\{(\Psi^{*})^{k}\}_{k=1}^{\infty}$ is a compact linear
skew-product dynamical system covering $\phi^{-1}$.

\proclaim{Theorem 1.4}
Equation (1a)+(1b) generates a compact linear
skew-product dynamical system $\{\Psi^{k}\}_{k=1}^{\infty}$
(\resp $\{\Pshat^{k}\}_{k=1}^{\infty}$) on the Banach
bundle $\CalB:=\Bb\times\Lone$ (\resp
$\CalBhat:=\Bb\times\Cone$), where $\phi(b):=b\cdot1$ for
$b\in\Bb$, and
$$
\align
&\psi(b)u_{0}:=u(1;b,u_{0})
\quad\text{for }(b,u_{0})\in\CalB.\\
(\text{\resp}&\pshat(b)u_{0}:=u(1;b,u_{0})
\quad\text{for }(b,u_{0})\in\CalBhat.)
\endalign
$$
The system $\{\Psi^{k}\}_{k=1}^{\infty}$ {\rm factorizes
through\/} $\{\Pshat^{k}\}_{k=1}^{\infty}$ in the sense
that $\psi(b)u\in\Cone$ for each $(b,u)\in\CalB$ and
$$
\pshat(b)=\psi(b)\circ i\quad\text{for each }b\in\Bb,
$$
where $i$ stands for the natural embedding
$\Cone\subset\Lone$.
\endproclaim
\demo{Proof}
The continuity of the assignment $b\mapsto\psi(b)$ or
$b\mapsto\pshat(b)$ as a mapping from $\Bb$ into the space
of bounded linear operators on $\Lone$ or $\Cone$ with the
strong operator topology follows by Theorem~1.3.  The
compactness follows by Proposition~1.2 and by the fact that
for $p$ sufficiently large the space $\Lp^{3/4}$ embeds
continuously in $C^{\mu}(\overline\Omega)$.  We derive the
continuity in the uniform operator topology along the lines
of Thm.~2.3.2. in Pazy [10].
\qed\enddemo

An immediate consequence is that for the dual skew-product
dynamical system, $\psi^{*}(b):\Cone^{*}\to
L^{\infty}(\Omega)$ for each $b\in\Bb$.  However, by the
Green's formula it follows that $\psi^{*}(b)u^{*}$,
$u^{*}\in\Cone^{*}$, equals the value at time $1$ of a
solution of the adjoint equation.  Since the adjoint
equation satisfies all the assumptions, we have the
following
\proclaim{Theorem 1.5}
The linear skew-product dynamical system
$\{(\Psi^{*})^{k}\}_{k=1}^{\infty}$ dual to the system
$\{\Psi^{k}\}$ generated by equation (1a)+\allowbreak(1b)
has the property that $\psi^{*}(b):\Cone^{*}\to\Cone$ for
each $b\in\Bb$.
\endproclaim

\head 2. Order and monotonicity \endhead

For a Banach space $X$ consisting of functions defined on
$\Omega$, by $X_{+}$ we denote the {\it cone\/} of
nonnegative functions. It is straightforward that $X_{+}$
is a closed convex set such that

a) For each $u\in X_{+}$ and $\alpha\ge0$ one has
${\alpha}u\in X_{+}$, and

b) $X_{+}$ contains no one-dimensional subspace.

The cones in $\Lp$, $1\le p\le\infty$, as well as in $\Cone$
are {\it generating\/}, that is, $X_{+}+X_{+}=X$.

We write $v\le u$ if $u-v\in X_{+}$, and $v<u$ if $v\le u$
and $v\ne u$.

A cone $X_{+}$ is called {\it solid\/} if its interior
$X_{++}$ is nonempty.  We write $v\ll u$ if $u-v\in
X_{++}$.  The standard cone in $\Lp$ is solid if and only
if $p=\infty$.  The standard cone in $\Cone$ is solid.

We say that a bounded linear functional $u^{*}\in X^{*}$ is
{\it nonnegative\/} if $\langle u^{*},u\rangle\ge0$ for
each $u\in X_{+}$, where $\langle\cdot,\cdot\rangle$
denotes the duality pairing.  When the cone $X_{+}$ is
generating, the set of all nonnegative functionals forms a
cone $(X^{*})_{+}$, called the {\it dual cone\/}.

A functional $u^{*}\in(X^{*})_{+}$ is called {\it uniformly
positive\/} if there is a constant $K>0$ such that $\langle
u^{*},u\rangle\ge K\|u\|$ for each nonzero $u\in X_{+}$,
where $\|\cdot\|$ stands for the norm in $X$.  Among the
spaces $\Lp$, only $\Lone$ admits uniformly positive
functionals: They are represented by functions in
$L^{\infty}(\Omega)$ whose values are a.e. larger than some
$K>0$ (in other words, elements of
$L^{\infty}(\Omega)_{++}$).  For more on cones the reader
is referred to Amann's paper [1].

\proclaim{Theorem 2.1}
Let $\{\Psi^{k}\}_{k=1}^{\infty}$ be the linear
skew-product dynamical system generated on $\CalB$ by
(1a)+\allowbreak(1b).  There exists an invariant
decomposition $\CalB=\CalBone\oplus\CalBtwo$, such that

{\rm (i)} $\codim\CalBone=1$, and for each $b\in\Bb$ the
fiber $\CalBone(b)$ is represented as the nullspace of a
uniformly positive functional $v^{*}(b)\in\Lone^{*}$.

{\rm (ii)} $\CalBtwo\setminus\CalZ\subset\Bb\times
(\Lone_{+}\cup-\Lone_{+})$, where $\CalZ$ denotes the null
section of $\CalB$.

{\rm (iii)} The mapping $\Psi|\CalBtwo$ is a bundle
automorphism.

{\rm (iv)} There are constants $D\ge1$ and $0<\lambda<1$
such that
$$
\frac{\|\psi^{(k)}(b)v_{1}\|_{\Lone}}
{\|\psi^{(k)}(b)v_{2}\|_{\Lone}}\le
D{\lambda}^{k}
\frac{\|v_{1}\|_{\Lone}}{\|v_{2}\|_{\Lone}}
$$
for each $(b,v_{1})\in\CalBone$,
$(b,v_{2})\in\CalBtwo\setminus\CalZ$, and each
$k\in\naturals$.
\endproclaim
\demo{Proof}
By Theorem~1.3, for each $b\in\Bb$, each $u_{0}\in\Lone$
and each $t>0$ we have $u(t;u_{0},b)\in\Cone$.  Assume that
a nonzero $u_{0}\in\Lone_{+}$.  As $u(t;u_{0},b)$ may fail
to be a classical solution, we cannot apply the standard
parabolic strong maximum principle as put forward \eg in
Protter and Weinberger [12].  However, since $u$ as well as
its spatial derivatives are continuous, we can establish
the parabolic strong maximum principle by reasoning along
the lines of Chapters~8 and~9 in Gilbarg and Trudinger [6],
which enables us to show after the pattern of Thm\.~4.1 in
Hirsch [9] or Cor\.~7.2.3 in Smith [14] that
$u(t;u_{0},b)\gg0$ for $t>0$.

Applying Thm\.~1 in \Polacik and \Terescak [11] to the
linear skew-product dynamical system
$\{\Pshat^{k}\}_{k=1}^{\infty}$, we obtain the existence of
two $\Pshat$-invariant subbundles $\CalS$ and $\CalT$ of
$\CalBhat$ such that:

(A1) $\CalBhat=\Bb\times\Cone$ is the direct sum of $\CalS$
and $\CalT$.

(A2) $\dim{\CalS}=1$, and $\CalS\setminus\CalZhat\subset
\Bb\times(\Cone_{++}\cup-\Cone_{++})$, where $\CalZhat$
stands for the null section of $\CalBhat$.

(A3) For each $b\in\Bb$ the fiber $\CalT(b)$ can be
(uniquely) represented by a normalized functional
$\hat{v}^{*}(b)\in\Cone^{*}_{+}$.  In~particular,
$\CalT\cap(\Bb\times\Cone_{+})=\CalZhat$.

(A4) $\Pshat|\CalS$ is a bundle automorphism.

(A5) There are constants $d\ge1$ and $0<\lambda<1$ such
that
$$
\frac{\|\pshat^{(k)}(b)v_{1}\|_{\Cone}}
{\|\pshat^{(k)}(b)v_{2}\|_{\Cone}}\le
d{\lambda}^{k}
\frac{\|v_{1}\|_{\Cone}}{\|v_{2}\|_{\Cone}}
$$
for each $(b,v_{1})\in\CalT$, $(b,v_{2})\in\CalS$,
$v_{2}\ne0$ and each $k\in\naturals$.

Now, the desired bundle $\CalBtwo$ equals simply
$\CalS$ considered a subbundle of the bundle
$\CalB=\Bb\times\Lone$.

In order to construct $\CalBone$, notice first that the
mapping $\Bb\ni b\mapsto\hat{v}^{*}(b)\in\Cone^{*}$ is
continuous.  Moreover, the $\Pshat$-invariance of $\CalT$
means that
$\pshat^{*}(\phi^{-1}(b))\hat{v}^{*}(\phi^{-1}(b))=
{\gamma}\hat{v}^{*}(b)$ for some $\gamma=\gamma(b)>0$. From
Theorem~1.5 we derive that for each $b\in\Bb$ the
functional $\hat{v}^{*}(b)$ is represented as a function
from $\Cone$ depending continuously on $b$.  Denote by
$v^{*}(b)$ the functional $\hat{v}^{*}(b)$ viewed as an
element of $\Cone$.  Applying the parabolic strong maximum
principle to the adjoint equation we get that
$v^{*}(b)\in\Cone_{++}$, consequently it is a uniformly
positive functional from $\Lone^{*}$.

To prove (iii) and (iv), use (A4), (A5) and the fact that
since $\CalBone$ has finite dimension, both the
$\Cone$-norm and the $\Lone$-norm on it are equivalent.
\qed\enddemo

The property described in (iv) is usually referred to as
{\it exponential separation\/} (continuous separation in
[11]).  The direct sum decomposition
$\CalB=\CalBone\oplus\CalBtwo$ uniquely defines the bundle
projection $P$ with image $\CalBtwo$ and kernel $\CalBone$.
The exponential separation can be formulated in the
following way:

There are constants $D\ge1$ and $0<\lambda<1$ such
that
$$
\frac{\|(\Id-P(\phi^{k}(b)))\psi^{(k)}(b)v\|_{\Lone}}
{\|P(\phi^{k}(b))\psi^{(k)}(b)v\|_{\Lone}}\le
D{\lambda}^{k}
\frac{\|(\Id-P(b))v\|_{\Lone}}{\|P(b)v\|_{\Lone}}
$$
for each $(b,v)\in\CalB\setminus\CalBone$ and each
$k\in\naturals$.

For $t\in\reals$ and $b\in\Bb$ put $\phi_{t}b$ to be
$b\cdot t$.  The family $\{\phi_{t}\}_{t\in\reals}$ forms a
flow (= continuous-time dynamical system) on the compact
metrizable space $\Bb$.  Define for $t\ge0$ and $b\in\Bb$ a
linear operator $\psi(t,b)\in{\Cal L}(\Lone)$ by the
formula
$$
\psi(t,b)u_{0}:=u(t;b,u_{0})\quad\text{for }u_{0}\in\Lone.
$$
For each $t\ge0$ the mapping $\Psi_{t}$ defined as
$$
\Psi_{t}(b,u):=(\phi_{t}b,\psi(t,b)u)
$$
is a bundle endomorphism of $\CalB$.  The family
$\{\Psi_{t}\}_{t\ge0}$ forms a {\it linear skew-product
semiflow\/} on the product bundle $\CalB$.  A consequence
of the semiflow axioms is the following {\it cocycle
identity\/}
$$
\psi(t_{1}+t_{2},b)=\psi(t_{1},\phi_{t_{2}}b)=
\psi(t_{2},\phi_{t_{1}}b)\quad\text{for }t_{1}\ge0,
t_{2}\ge0. \tag{2.1}
$$
For basic properties of linear skew-products semiflows the
reader is referred to Sacker and Sell [13] or to Chow and
Leiva [3].

\proclaim{Theorem 2.2}
Let $\{\Psi_{t}\}_{t\ge0}$ be the linear skew-product
semiflow generated on the product bundle
$\CalB=\Bb\times\Lone$ by Equation (1a)+\allowbreak(1b).
Then

{\rm (i)} The subbundles $\CalBone$ and $\CalBtwo$ are
$\Psi_{t}$-invariant in the sense that if $(b,u)\in\CalBi$
then $(\phi_{t}b,\psi(t,b)u)\in\CalBi$ for $t\ge0$, where
$i=1$, $2$.

{\rm (ii)} $\{\Psi_{t}|\CalBtwo\}$ extends uniquely to a
linear skew-product flow on $\CalBtwo$.

{\rm (iii)} There are constants $D'\ge1$
and $\mu>0$ ($\mu=-\log\lambda$) such
that
$$
\frac{\|\psi(t,b)v_{1}\|_{\Lone}}
{\|\psi(t,b)v_{2}\|_{\Lone}}\le
D'e^{-{\mu}t}
\frac{\|v_{1}\|_{\Lone}}{\|v_{2}\|_{\Lone}}
$$
for each $(b,v_{1})\in\CalBone$,
$(b,v_{2})\in\CalBtwo\setminus\CalZ$, and each
$t\ge0$.
\endproclaim
\demo{Proof}
We start by proving (iii).  Indeed, by the cocycle identity
we have $\psi(t,b)=\psi(t-[t],\phi(t,b))\circ\psi([t],b)$,
where $[t]$ stands for the integer part of $t$.  By the
standard argument (compare \eg the proof of Lemma~3.3 in
Sacker and Sell [13]) there is a positive constant $C'$
such that $\|\psi(s,b)u\|_{\Lone}\le C'\|u\|_{\Lone}$ for
$s\in[0,1]$, $b\in\Bb$ and $u\in\Lone$.  Further, as
$\CalBtwo$ is one-dimensional, there is a positive constant
$C''$ such that $\|\psi(s,b)u\|_{\Lone}\ge
C''\|u\|_{\Lone}$ for $s\in[0,1]$ and $(b,u)\in\CalBtwo$.
From this it follows that for $(b,v_{1})\in\CalBone$,
$(b,v_{2})\in\CalBtwo\setminus\CalZ$, $t\ge0$, we have
$$
\frac{\|\psi(t,b)v_{1}\|_{\Lone}}
{\|\psi(t,b)v_{2}\|_{\Lone}}\le
\frac{C'\|\psi([t],b)v_{1}\|_{\Lone}}
{C''\|\psi([t],b)v_{2}\|_{\Lone}}\le
\frac{C'D}{C''}{\lambda}^{[t]}
\frac{\|v_{1}\|_{\Lone}}{\|v_{2}\|_{\Lone}}\le
\frac{C'D}{C''}{\lambda}^{t-1}
\frac{\|v_{1}\|_{\Lone}}{\|v_{2}\|_{\Lone}}.
$$
Putting $D':=C'C/{\lambda}C''$ and $\mu:=-\log\lambda$ gives
the desired result.

Part (ii) is a consequence of the fact that $\Psi|\CalBtwo$
is a bundle automorphism and the cocycle property.

The proof of the $\Psi_{t}$-invariance of $\CalBtwo$ is
straightforward:  $(b,u)\in\CalBtwo$ is equivalent to
$(\Id-P(b))u=0$, which yields, by part (iii),
$(\Id-P(\phi_{t}b))\psi(t,b)u$ for all $t\ge0$.  In order
to establish the $\Psi_{t}$-invariance of $\CalBone$,
suppose by way of contradiction that for some
$(b',v)\in\CalBone$ and some $t'>0$, $t'\notin\naturals$,
we have $\psi(t',b')v\notin\CalBone$.  This means that
$P(\phi_{t'}b')v\ne0$.  As
$\psi([t']+1,b')v=\psi([t']+1-t',\phi_{t'}b')v$, from (iii)
we deduce that the ratio
$$
\frac{\|(\Id-P(\phi_{[t']+1}b'))
\psi([t']+1,b')v\|_{\Lone}}
{\|P(\phi_{[t']+1}b')\psi([t']+1,b')v\|_{\Lone}}
$$
is finite, which contradicts the fact that
$(\phi_{[t']+1}b',\psi([t']+1,b')v)\in\CalBone$.
\qed\enddemo

\proclaim{Theorem 2.3}
Let $u:\reals\times\overline\Omega\to\reals$ be a nonzero
global solution to (1a)+\allowbreak(1b) such that
$u(t,x)\ge0$ for each $t\in\reals$ and each
$x\in\overline\Omega$.  Then for each $t\in\reals$ the pair
$(a_{0}\cdot t,u(t,\cdot))$ belongs to the one-dimensional
subbundle $\CalBtwo$.
\endproclaim
\demo{Proof}
An application of the parabolic strong maximum principle
yields $u(t,\cdot)\in\Cone_{++}$ for each $t\in\reals$.

We claim that there is a positive constant $L$ such that
$$
\|P(b)w\|_{\Lone}\ge L\|w\|_{\Lone}\quad
\text{for each }b\in\Bb\text{ and each }w\in\Lone_{+}.
\tag{2.2}
$$
Indeed, as the fiber $\CalBone(b)$ is the nullspace of the
functional $\hat{v}^{*}(b)$, one has
$\langle\hat{v}^{*}(b),w\rangle=
\langle\hat{v}^{*}(b),P(b)w\rangle$ for all $w\in\Lone$.
Further, since $\CalBtwo$ has dimension one and the
functional $\hat{v}^{*}(b)$ is uniformly positive, there is
a positive number $h(b)$ such that
$\|w\|_{\Lone}=h(b)\langle\hat{v}^{*}(b),w\rangle$ for all
$w\in\Lone_{+}$ with $(b,w)\in\CalBtwo$.
Consequently $\|P(b)w\|_{\Lone}=
h(b)\langle\hat{v}^{*}(b),P(b)w\rangle$ for all
$w\in\Lone_{+}$.  The positive function $h(\cdot)$ is
easily seen to be continuous, so there is $h>0$,
$h=\min\{h(b):b\in\Bb\}$, such that
$\|P(b)w\|_{\Lone}\ge h\langle\hat{v}^{*}(b),w\rangle$ for
all $b\in\Bb$ and $w\in\Bb$.  Now it remains to notice that
by the continuity of the mapping $\Bb\ni
b\mapsto\hat{v}^{*}(b)\in\Lone^{*}$ and the fact that
uniformly positive functionals form an open set in
$\Lone^{*}$, the positive constant
in the definition of uniform positivity can be chosen
independent of $b$.  Formula (2.2) follows immediately.

As a consequence, for each nonzero $w\in\Lone_{+}$ and each
$b\in\Bb$ one has
$$
\frac{\|(\Id-P(b))w\|_{\Lone}}{\|P(b)w\|_{\Lone}}\le
\frac{1+N}{L},
$$
where $N:=\sup\{\|P(b))w\|_{\Lone}:b\in\Bb,
\|w\|_{\Lone}=1\}$.

For each $t\in\reals$ denote $\tilde{u}(t)=u(t,\cdot)$
regarded as an element of $\Lone$.  Suppose to the contrary
that there is $t'$ such that $(a_{0}\cdot t',\tilde{u}(t'))$
does not belong to $\CalBtwo$.  Put
$$
M:=\frac{\|(\Id-P(a_{0}\cdot t'))\tilde{u}(t')\|_{\Lone}}
{\|P(a_{0}\cdot t')\tilde{u}(t')\|_{\Lone}}>0,
$$
and
$$
t'':=t'+\frac{1}{\mu}\log{\frac{ML}{2D'(1+N)}},
$$
where $\mu$ and $D'$ are constants from Theorem~2.2(iii).
As $M\le(1+N)/L$, we have $t''<t'$.

An application of Theorem~2.2(iii) yields
$$
\align
&\quad\frac{\|(\Id-P(a_{0}\cdot t'))\tilde{u}(t')\|_{\Lone}}
{\|P(a_{0}\cdot t')\tilde{u}(t')\|_{\Lone}}\\
&=
\frac{\|\psi(t'-t'',a_{0}\cdot t'')(\Id-P(a_{0}\cdot t''))
\tilde{u}(t'')\|_{\Lone}}
{\|\psi(t'-t'',a_{0}\cdot t'')P(a_{0}\cdot t'')
\tilde{u}(t'')\|_{\Lone}}\\
&\le
D'e^{-{\mu}(t'-t'')}
\frac{\|\psi(t'-t'',a_{0}\cdot t'')(\Id-P(a_{0}\cdot t''))
\tilde{u}(t'')\|_{\Lone}}
{\|\psi(t'-t'',a_{0}\cdot t'')P(a_{0}\cdot t'')
\tilde{u}(t'')\|_{\Lone}}\\
&\le
D'e^{-{\mu}(t'-t'')}\frac{1+N}{L}\\
&\le\frac{M}{2},
\endalign
$$
a contradiction.
\qed\enddemo

Finally, we formulate now our main result.
\proclaim{Corollary 2.4}
Assume that $u_{1}$, $u_{2}$ are nonzero global solutions
of (1a)+\allowbreak(1b) such that $u_{1}(t,x)\ge0$ and
$u_{2}(t,x)\ge0$ for all $t\in\reals$ and all
$x\in\overline\Omega$.  Then there is a positive constant
$\kappa$ such that $u_{1}(t,x)={\kappa}u_{2}(t,x)$ for all
$t\in\reals$ and all $x\in\overline\Omega$.
\endproclaim

\enddocument